\documentclass{article}
\usepackage{maa-monthly}

\usepackage{verbatim}
\usepackage{amsmath}

\theoremstyle{theorem}
\newtheorem{theorem}{Theorem}

\theoremstyle{definition}

\makeatletter
\def\cosec{\mathop{\operator@font cosec}\nolimits}

\begin{document}

\title{Geometric Approach to the integral $\int \sec x\,dx$}
\markright{Geometric Approach to $\int \sec x\,dx$}
\author{Udita Katugampola}

\maketitle

\begin{abstract}
We give a geometric proof of the evaluation of the integral $\int \sec x\,dx$ which is normally done using a rather ad hoc approach. 
\end{abstract}

\section{Introduction}
As known to the mathematics community, the evaluation of two familiar integrals, namely, $\int \sec x\, dx$ and $\int \cosec x\, dx$ follow a rather ad hoc approach. The catch is that the integrand is rewritten in a manner that the transformed integral takes the form $\int\frac{u^\prime}{u}\, dx$ and then use the fact that $\int \frac{1}{u} du = \ln |u| + C$. For example, in the case of $\int \sec x\, dx$ the numerator and the denominator of the integrand are multiplied 
by a factor of $ \sec x + \tan x $ which brings the integral into the form in question. 
The draw back of the proceeding method is that it kind of require the knowledge of the solution even before it is evaluated. To over come this difficulty, Jacob and Osler \cite{osler1} proposed an approach which use a very interesting geometric result given below.

\begin{theorem}\label{thm1}
The area of triangle $OAB$ in Figure 1 equals the area of trapezoid $ABDC$.
\end{theorem} 

\begin{figure}[htbp]
	\centering
		\includegraphics[height=2.4in, width=2.4in]{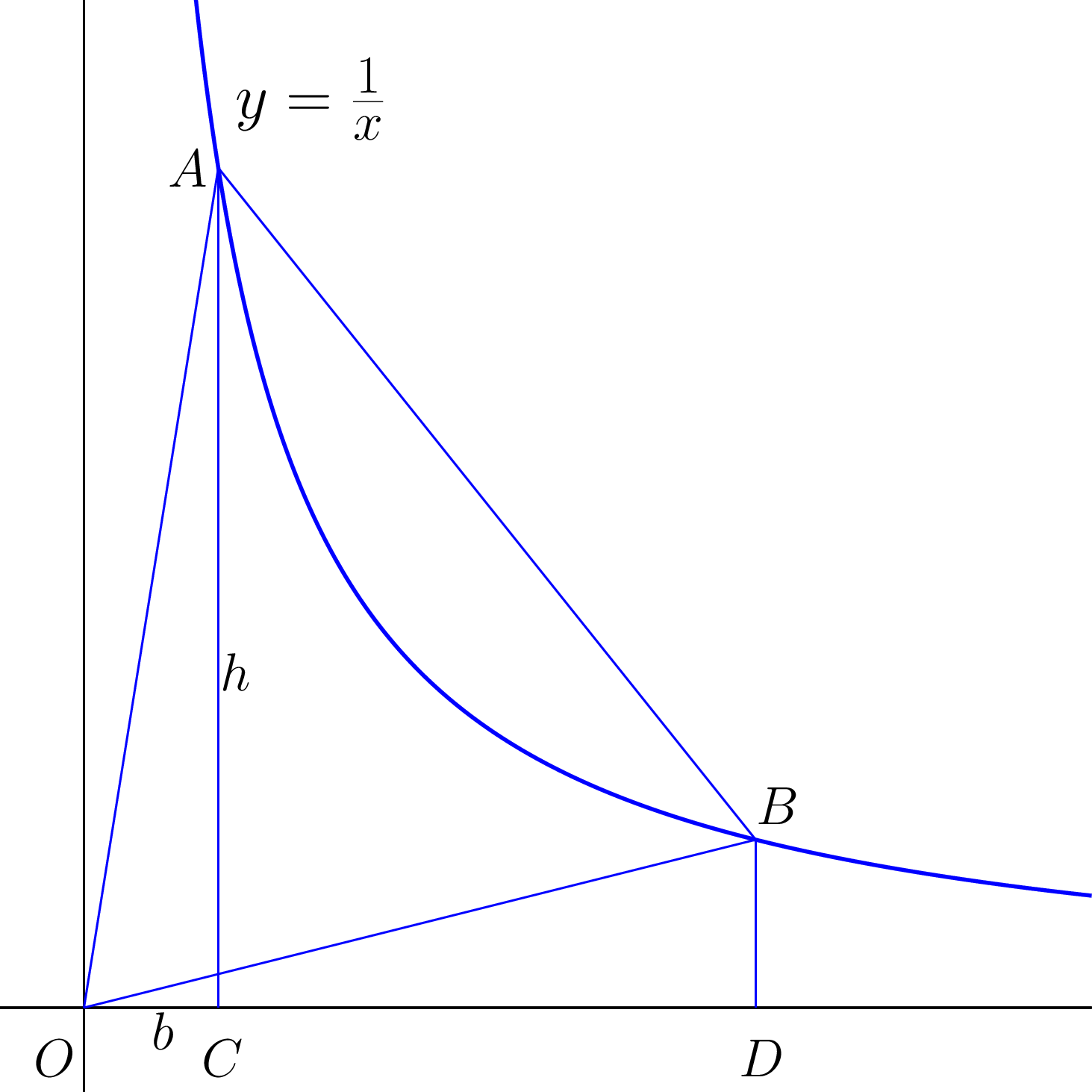}
	\caption{Equal Areas}
	\label{fig:Figure-1}
\end{figure}

This interesting result is not hard to prove by observing that the product of the height ($h$) and base ($b$) of any triangle such as $AOC$ or $BOD$ is a constant, thus making their areas keep fixed at $\frac{1}{2}bh = \frac{1}{2}$, since one vertex of such a triangle is always lies on the curve $y=\frac{1}{x}$. This observation is first made by James Gregory in 1667, according to Osler \cite{osler1}. 

Theorem~\ref{thm1} also suggests the following very important result \cite{osler1}.

\begin{theorem}\label{thm2}
The area of sector $OAB$ under the curve $y=\frac{1}{x}$ is equal to the area $ABDC$ under the same curve.
\end{theorem} 

The authors in \cite{osler1} use this observation to provide a novel geometric proof to the following identity which is normally proved using a similar ad hoc approach. 

\begin{theorem}\label{thm3}
$\int \csc x\,dx = -\ln(\csc x + \cot x) +C.$
\end{theorem} 

 This raises the question that whether we can proceed with a similar approach to the twin identity for $\sec x$. The purpose of the rest of this paper is to extend this idea to provide a geometric approach to the following sibling identity.

\begin{theorem}\label{thm4}
$\int \sec x\,dx = \ln(\sec x + \tan x) +C.$
\end{theorem} 

The geometric proof of this result requires some additional observations. First we notice that as the integral $\int \csc x\,dx$ corresponds to the curve $y=\frac{1}{x}$, the integral $\int \sec x\,dx$ corresponds to that of $x^2 - y^2 = 2$. Thus the method of the proof consists of first transform the coordinates and then use Theorem~\ref{thm2} as in the case of \cite{osler1}. 

 \section{Coordinate Transformation}
First consider the familiar coordinate transformation which rotates the $xy-$Cartesian plane about the origin through an angle of $\theta$ in the counter-clockwise direction to obtain the new axis of coordinates $X-Y$. 
\begin{figure}[htbp]
	\centering
		\includegraphics[height=2.4in, width=2.4in]{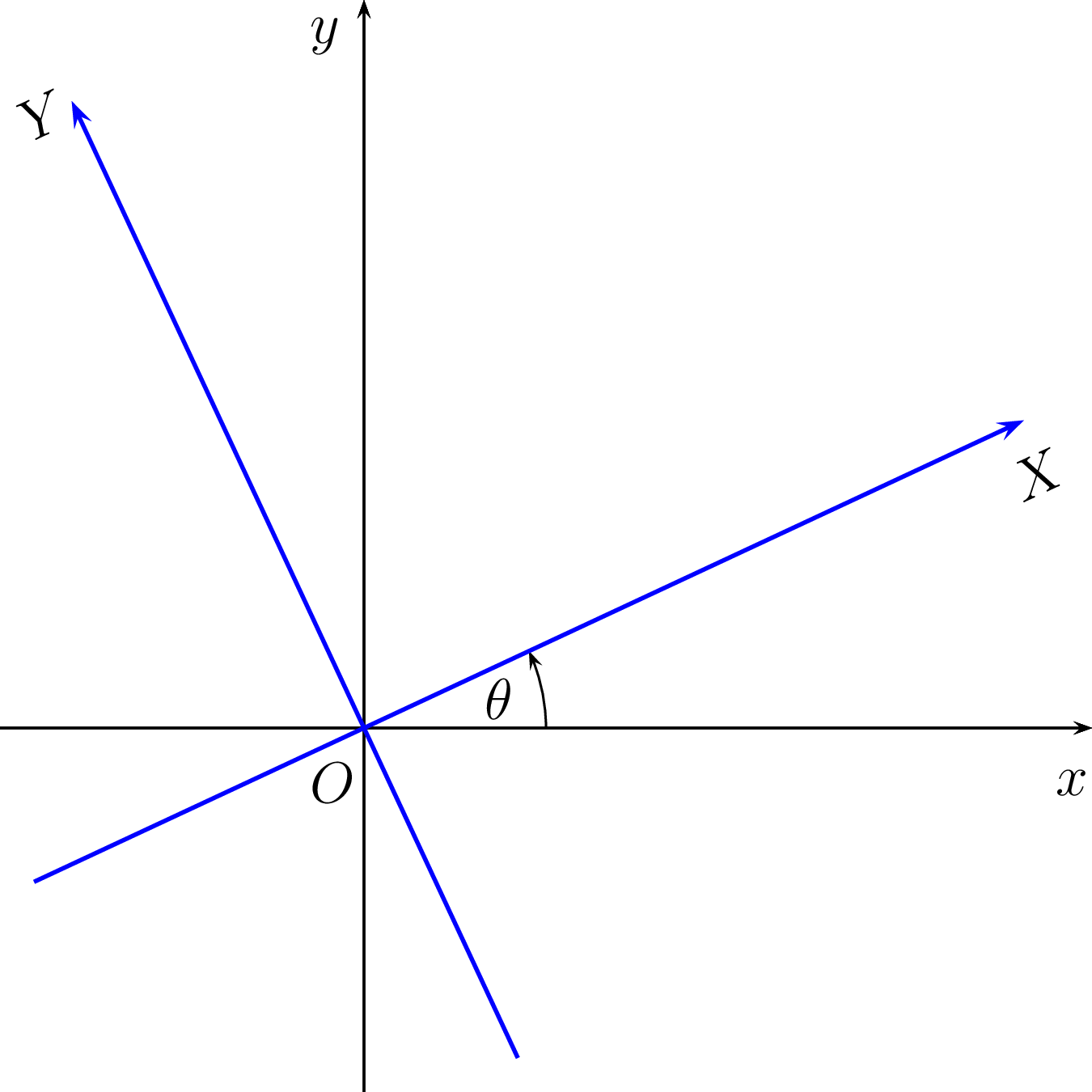}
	\caption{Axes Rotation}
	\label{fig:Figure-2}
\end{figure}

 The relation between the two coordinate systems is given by
\begin{equation}\label{tran1}
   \begin{bmatrix}
             X \\
             Y 
					\end{bmatrix} 
		 = \left(\begin{array}{cc}\cos \theta & \sin \theta \\
		                                 -\sin \theta &  \cos \theta 
										\end{array} \right) \begin{bmatrix}
             x \\
             y 
		\end{bmatrix} 
\end{equation}

For example, consider the rotation of $xy-$coordinates about the origin through an angle of $\theta = \frac{\pi}{4}$ in the clockwise direction. According to the Eq. (\ref{tran1}), the transformed coordinates are given by

\begin{align}
   X &= \frac{1}{\sqrt{2}}(x-y) \\
	 Y &= \frac{1}{\sqrt{2}}(x+y).
\end{align}

This, in turn, transforms the curve $x^2 -y^2 =2$ in the $xy-$plane to the curve given by
\vspace{-0.3cm}
\begin{equation}
   XY = \frac{1}{\sqrt{2}}(x-y)\cdot \frac{1}{\sqrt{2}}(x+y) = 1
\end{equation}
in the $XY-$plane. The equation of the new curve takes the form $XY=constant$, thus allowing us to use the Theorem~\ref{thm2} in the transformed coordinate system. This result would be the key ingredient in proving Theorem~\ref{thm4} above.  

\begin{proof}[$\textbf{Proof of Theorem 4}$]
In Figure 3, let $\angle BOE = \alpha$ and $\angle AOE = \beta$. Then the area of the sector $AOB$, using polar coordinates, is $\int_\alpha^\beta \frac{r^2}{2} d\vartheta $. But since the points $A$ and $B$
\begin{figure}[h]
	\centering
		\includegraphics[height=2.5in, width=2.8in]{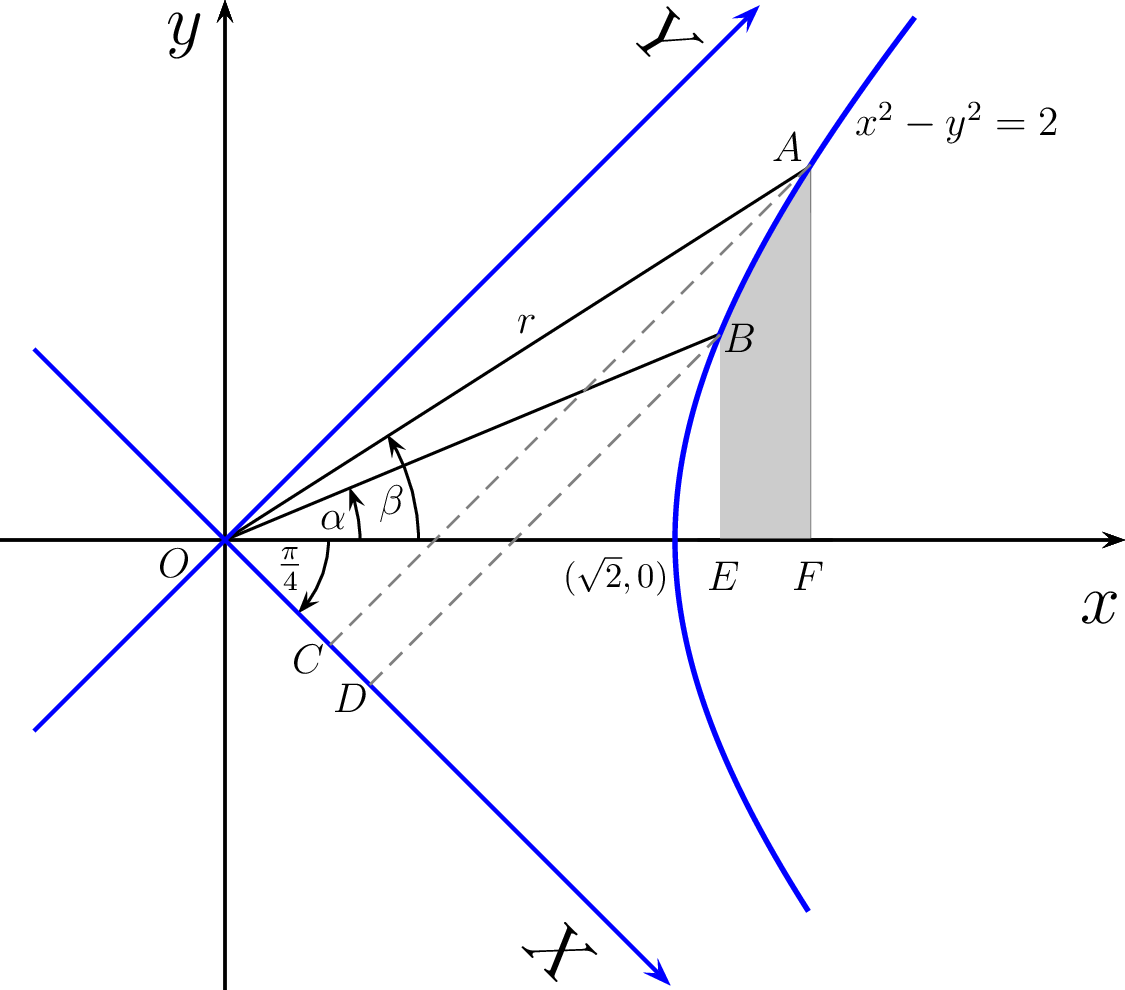}
	\caption{Evaluating the integral $\int \sec x dx$}
	\label{fig:Figure-3}
\end{figure}
\hspace{-.3cm}
\noindent
~lie on the curve $x^2 - y^2 = 2$, we have $r^2\cos^2 \vartheta - r^2 \sin^2 \vartheta =2$, so that $r^2=\frac{2}{\cos^2 \vartheta - \sin^ 2 \vartheta} = 2\sec 2\vartheta$. Thus, 
\hspace{-.5cm}
\[
    \mbox{Area of sector } AOB =  \int_\alpha^\beta \sec 2\vartheta d\vartheta.
\]
But, by Theorem~\ref{thm2}, this must be equal to the area of the region $ABCD$ which we can evaluate in the transformed coordinate system. Without loss of generality, let the point $B$ be $(\sqrt{2},0)$ in the $xy-$coordinate system so that $X_D=1$. In the $XY-$coordinate system the area in question is equal to $\int_{X_C}^{X_D} \frac{1}{X}dX = \ln\frac{X_D}{X_C} = -\ln X_C$. But notice that
\begin{align*}
   X_C = r\cos (\beta + \frac{\pi}{4}) &= \frac{\cos (\beta + \frac{\pi}{4})}{\sqrt{\sin (\beta + \frac{\pi}{4})\cos (\beta + \frac{\pi}{4})}} 
	                                     =\frac{\sqrt{\cos^2 (\beta + \frac{\pi}{4})}}{\sqrt{\sin (\beta + \frac{\pi}{4})\cos (\beta + \frac{\pi}{4})}}\\
																			 &=\sqrt{\frac{1+\cos (2\beta + \frac{\pi}{2})}{\sin (2\beta + \frac{\pi}{2})}} 
																			= \sqrt{\frac{1 - \sin 2\beta}{\cos 2\beta}}\\
																			&= \sqrt{\sec 2\beta - \tan 2\beta}.
\end{align*}
Thus, we have 
\[
    \int_{0}^\beta \sec 2\vartheta d\vartheta = -\ln X_C = -\ln \sqrt{\sec 2\beta - \tan 2\beta} = \frac{1}{2}\ln (\sec 2\beta + \tan 2\beta).
\]
Now letting $ z = 2\vartheta$, we have 
\[
    \int_{0}^x \sec z\,dz = \ln (\sec x + \tan x),
\]
as required. This completes the proof of the main result. 
\end{proof}


\begin{acknowledgment}{Acknowledgment}
Research of the author was supported by the Department of Defense (USA) grant 67459MA-15-139 MJ.
\end{acknowledgment}

\begin{affil}
Department of Mathematical Sciences, University of Delaware, Newark DE 19716\\
uditanalin@yahoo.com
\end{affil}

\vfill\eject

\end{document}